\documentclass[12pt,reqno,a4wide]{article}
\usepackage{makeidx}
\usepackage{comment}
\usepackage{mathptmx}
\usepackage{a4wide}
\usepackage{amssymb}
\usepackage{amsfonts}
\usepackage{amsthm}
\usepackage{amsmath}
\usepackage{dsfont}
\usepackage{graphicx}
\usepackage{shadow}
\usepackage[all]{xy}
\usepackage{enumerate}
\usepackage{mathrsfs}
\usepackage{upgreek}
\usepackage{stmaryrd}
\usepackage[utf8]{inputenc}
\usepackage{array}
\usepackage{tikz-cd}
\usepackage{thumbpdf}
\usepackage[colorlinks=true,pagebackref]{hyperref}
\usepackage{textcomp}
\usepackage{chapterbib}
\usepackage{graphics}
\usepackage{longtable}
\usepackage{epsf}
\usepackage{bm}
\usepackage{adjustbox}

\newtheorem{theorem}{Theorem}[section]
\newtheorem{lem}[theorem]{Lemma}
\newtheorem{thm}[theorem]{Theorem}
\newtheorem{prop}[theorem]{Proposition}
\newtheorem{cor}[theorem]{Corollary}
\theoremstyle{definition}
\newtheorem*{Beweis}{Proof}
\newtheorem{defn}[theorem]{Definition}
\newtheorem{ex}[theorem]{Example}
\newtheorem{exs}[theorem]{Examples}
\newtheorem{defns}[theorem]{Definitions}
\newtheorem{rem}[theorem]{Remark}
\newtheorem{rems}[theorem]{Remarks}
\newtheorem{punto}[theorem]{}

\input{tcilatex}
\swapnumbers
\CompileMatrices
\input{tcilatex}
\begin{document}

\title{PS-Hollow Representations of Modules over Commutative Rings
}
\author{ $%
\begin{array}{ccc}
\text{Jawad Abuhlail}\thanks{Corresponding Author; Email: abuhlail@kfupm.edu.sa. \newline
The authors would like to acknowledge the support provided by the Deanship
of Scientific Research (DSR) at King Fahd University of Petroleum $\&$
Minerals (KFUPM) for funding this work through projects No. RG1213-1 $\&$
RG1213-2} &  & \text{Hamza Hroub}\thanks{%
The paper is extracted from the Ph.D. dissertation of Dr. Hamza Hroub under
the supervision of Prof. Jawad Abuhlail} \\
\text{Department of Mathematics and Statistics} &  & \text{Department of
Mathematics } \\
\text{King Fahd University of Petroleum $\&$ Minerals}\&\text{ Minerals} &  & \text{Arab American University } \\
\text{31261 Dhahran, KSA} &  & \text{Jenin, Palestine}%
\end{array}%
$}
\date{\today }
\maketitle

\begin{abstract}
Let $R$ be a commutative ring and $M$ a non-zero $R$-module. We introduce the class of \emph{pseudo strongly hollow submodules} (\emph{PS-hollow submodules}, for short) of $M$. Inspired by the theory of modules with \emph{secondary representations}, we investigate modules which can be written as \emph{finite} sums of PS-hollow submodules. In particular, we provide existence and uniqueness theorems for the existence of \emph{minimal} PS-hollow strongly representations of modules over Artinian rings.
\end{abstract}

\section*{Introduction}

This paper is part of our continuing project of investigating the different notions of primeness and coprimeness for (sub)modules of a given a non-zero module $M$ over a (commutative) ring $R$ in their \emph{natural context} as prime (coprime) elements in the lattice $Sub_R(M)$ of $R$-submodules with the canonical action of the poset $Ideal(R)$ of ideals of $R$. This approach proved to be very appropriate and enabled use to prove several results in this general setting and to provide more elegant and shorter proofs of our results. Moreover, it enabled us to generalize several notions and dualize them in a more systematic and elegant way.

Generalizing the notion of a \emph{strongly hollow element} in a lattice, we introduce for a lattice with an \emph{action of a poset} the notion of a \emph{pseudo strongly hollow element}. The two notions are equivalent in case the lattice is \emph{multiplication}. Considering the lattice $Sub_R(M)$ of a non-zero module $M$ over a commutative ring $R$, we obtain new class of modules, which we call \emph{pseudo strongly hollow modules}. We study this class of $R$-modules, as well as modules which can be written as \emph{finite sums} of their pseudo strongly hollow submodules. In particular, we provide existence and uniqueness theorems of such representation over Artinian rings.

This paper consists of two sections. In Section 1, we define, for a bounded lattice $\mathcal{L} = (L,\wedge, \vee, 0, 1)$, several notions of primeness for elements in $L\backslash \{1\}$ as well as several coprimeness notions for elements in $L\backslash \{0\}$. In Theorem \ref{Proposition 5.8}, we prove that the spectrum $Spec^{c}(\mathcal{L})$ of coprime elements in $L$ is nothing but the spectrum $Spec^{s}(\mathcal{L}^{0})$ of second elements in the dual bounded lattice $\mathcal{L}^{0}:=(L,\vee ,\wedge ,1,0).$

In Section 2, we apply the results of Section 1 to the lattice $\mathcal{L}:=Sub_R(M)$ of
submodules of a non-zero module $M$ over a commutative ring $R.$ We present the
notion of a \emph{pseudo strongly hollow submodule }(\emph{PS-hollow submodule} for
short) $N\leq M$ as dual to the \emph{pseudo strongly irreducible submodules}. Modules which are
finite sums of PS-hollow submodules are said to be \emph{PS-hollow\ representable}.
Proposition \ref{Remark 5.6} asserts the existence of \emph{minimal} PS-hollow
representations for PS-hollow representable modules over Artinian rings. The
First and the Second Uniqueness Theorems of minimal pseudo strongly hollow
representations are provided in Theorems \ref{Theorem 5.9} and \ref{Theorem
5.10}, respectively. Sufficient conditions for $_{R}M$ to have a PS-hollow
representation are given in Proposition \ref{hollow-PS}. Finally, Theorem %
\ref{Theorem 5.17} investigates semisimple modules each PS-hollow submodules
of which is simple.

\section{Primeness and Coprimeness Conditions for Lattices}

\hspace{1cm} In this section, we provide some preliminaries and study several notions of \emph{primeness} and \emph{coprimeness} for elements in a complete lattice $\mathcal{L}:=(L,\wedge ,\vee,0,1)$ attaining an action of a poset $(S,\leq ).$

Throughout, $S=(S,\leq)$ is a non-empty poset and $S^0=(S,\geq)$ is the dual poset.

\begin{punto}
(\cite{G}) A \emph{lattice}
\index{lattice} $\mathcal{L}$ is a \emph{poset} $(L,\leq )$ closed under two
binary commutative, associative and idempotent operations: $\wedge $ (\emph{%
meet}) and $\vee $ (\emph{join}), and we write $\mathcal{L}=(L,\wedge ,\vee)$; we say that  $\mathcal{L}$ is a \emph{bounded lattice} iff there exist $0,1 \in L$ such that $0 \leq x \leq 1$ for all $x\in L$.
We say that a lattice $(L,\wedge ,\vee )$ is a \emph{complete lattice} iff
$\bigwedge\limits_{x\in H}x$ and $\bigvee\limits_{x\in H}x$ exist in $L$ for any
$H\subseteq L$. Every complete lattice is bounded with $0= \bigwedge\limits_{x\in L}x$ and $1= \bigvee\limits_{x\in L}x$.

For two (complete) lattices $\mathcal{L}=(L,\wedge ,\vee )$ and
$\mathcal{L}^{\prime }=(L^{\prime },\wedge ^{\prime },\vee ^{\prime }),$ a
\emph{homomorphism of (complete) lattices}
\index{homomorphism of lattices} from $\mathcal{L}$ to $\mathcal{L}^{\prime
} $ is a map $\varphi :L\longrightarrow L^{\prime }$ that preserves finite
(arbitrary) meets and finite (arbitrary) joins.
\end{punto}

The notion of a \emph{strongly hollow submodule} was introduced by Abuhlail in \cite{Abu2015}, as dual to that of \emph{strongly irreducible submodules}. The notion was generalized to general lattices and investigated by Abuhlail and Lomp in \cite{AbuC0}.

\begin{punto}
Let $\mathcal{L}=(L,\wedge ,\vee ,0,1)$ be a bounded lattice.

\begin{enumerate}
\item An element $x\in L\backslash \{1\}$ is said to be:

\emph{irreducible} (or \emph{uniform}) iff for any $a,b\in L$ with $a\wedge b=x,$
we have $a=x$ or $b=x$;

\emph{strongly irreducible} iff for any $a,b\in L$ with $%
a\wedge b\leq x,$ we have $a\leq x$ or $b\leq x$.

\item An element $x\in L\backslash \{0\}$ is said to be:

\emph{hollow }iff whenever for any $a,b\in L$ with $x=a\vee b,$ we have $x=a$
or $x=b;$

\emph{strongly hollow} iff for any $a,b\in L$ with $x\leq a\vee b,$ we have $x\leq a$ or $x\leq b$.
\end{enumerate}

We denote the set of irreducible (resp. strongly irreducible, hollow, strongly hollow) elements in $L$ by $I(\mathcal{L})$ (resp. $SI(\mathcal{L})$, $H(\mathcal{L})$, $SH(\mathcal{L})$).

We say that $\mathcal{L}$ is a \emph{hollow lattice} (resp. \emph{uniform lattice}) iff $1$ is hollow (0 is uniform).
\end{punto}

\begin{punto}
\index{lattice}
\label{action}Let $\mathcal{L}=(L,\wedge,\vee )$ be a lattice. An $S$-\textit{action}
\index{$S-$action} on $\mathcal{L}$ is a map $\rightharpoonup :S\times L\longrightarrow L$
satisfying the following conditions for all $s,s_{1},s_{2}\in S$ and $x,y\in L$:
\begin{enumerate}
  \item $s_{1}\leq _{S}s_{2}\Rightarrow s_{1}\rightharpoonup x\leq
s_{2}\rightharpoonup x$.
  \item $x\leq y\Rightarrow s\rightharpoonup x\leq s\rightharpoonup y$.
  \item $s\rightharpoonup x \leq x$.
\end{enumerate}
A bounded lattice $\mathcal{L}=(L,\wedge,\vee,0,1)$ with an $S$-action is \emph{multiplication}
iff for every element $x\in L$, there is some $s\in S$ such that $x=s\rightharpoonup 1$.
\end{punto}

\begin{ex}
\index{complete lattice}
\label{LAT}Let $M$ be an $R$-module. The complete lattice $LAT(_{R}M)$ of $R$-submodule has an $Ideal(R)$-action defined by the canonical product $IN$ of an ideal $I\leq R$ and a submodule $N\leq M$.
\end{ex}

\begin{rem}
\label{dual action}Let $\mathcal{L}=(L,\wedge,\vee ,0,1)$ a bounded lattice with an $S$-action $\rightharpoonup
: S\times L\longrightarrow L.$ The dual lattice $\mathcal{L}^{0}$  has an $S^0$-action given by
\begin{equation*}
s\rightharpoonup ^{0}x=(s\rightharpoonup 1)\vee x%
\text{ }  \label{s0x}, \text{ for all } s\in S \text{ and }x\in L.
\end{equation*}
\end{rem}

We generalized the notion of a \emph{strongly hollow element} of a lattice investigated by Abuhlail and Lomp in $\cite{AbuC0}$ to a \emph{strongly hollow element} of a lattice with an action from a poset. Moreover, we introduced its dual notion of a \emph{pseudo strongly irreducible element} which is a generalization of the notion of a \emph{strongly irreducible element}.

\begin{defns}
\label{spectra}Let $(\mathcal{L},\rightharpoonup )$ a bounded lattice with an $S$-action. We say that:

\begin{enumerate}
\item $x\in L\backslash \{1\}$ is

\emph{pseudo strongly irreducible} iff for all $y\in L$ and $s\in S:$
\begin{equation}
(s\rightharpoonup 1)\wedge y\leq x \hspace{8pt} \Rightarrow \hspace{8pt} s\rightharpoonup 1\leq x \text{\hspace{8pt}
or \hspace{8pt} }y\leq x;  \label{psi}
\end{equation}

\emph{prime} iff for all $y\in L$ and $%
s\in S$ with
\begin{equation}
s\rightharpoonup y\leq x\hspace{8pt} \Rightarrow \hspace{8pt} s\rightharpoonup 1\leq x\text{\hspace{8pt}
or \hspace{8pt} }y\leq
x.  \label{p}
\end{equation}

\emph{coprime} iff for all $s\in S:$%
\begin{equation}
s\rightharpoonup 1\leq x\text{\hspace{8pt}
or \hspace{8pt} }(s\rightharpoonup 1)\vee x=1  \label{c}
\end{equation}

\item $x\in L\backslash \{0\}$ is

\emph{pseudo strongly hollow} (or \emph{PS-hollow} for short) iff for all $s\in S:$
\begin{equation}
z\leq s\rightharpoonup x + y \Rightarrow z\leq s \rightharpoonup 1 \text{ or }z\leq y.  \label{psh}
\end{equation}%

\emph{second} iff for all $s\in S:$%
\begin{equation}
s\rightharpoonup x=x\text{\hspace{8pt}
or \hspace{8pt} }s\rightharpoonup x=0  \label{s}
\end{equation}

\emph{first} iff for all $y\in L$ and $s\in S$ with%
\begin{equation}
s\rightharpoonup y=0\text{\hspace{8pt}
and \hspace{8pt} }y\leq x\hspace{8pt} \Rightarrow \hspace{8pt} s\rightharpoonup x=0\text{\hspace{8pt}
or \hspace{8pt} }y=0.  \label{f}
\end{equation}%
\end{enumerate}

The spectrum of pseudo strongly irreducible (resp. prime, coprime, pseudo strongly hollow, second, first) elements of $\mathcal{L}$ is denoted by $Spec^{psi}(\mathcal{L})$ (resp. $Spec^{p}(\mathcal{L})$, $Spec^{c}(\mathcal{L})$, $Spec^{s}(\mathcal{L})$, $Spec^{f}(\mathcal{L})$).
\end{defns}

\begin{lem}
\index{complete lattice}
\label{00=*}Let $\mathcal{L}=(L,\wedge ,\vee,0,1)$ be a bounded lattice with an $S$-action and define%
\begin{equation}
s\rightharpoonup ^{\ast }x=(s\rightharpoonup 1)\wedge x  \label{s*x}
\end{equation}
for all $s\in S$ and $x\in L.$ Then $((\mathcal{L},\rightharpoonup
)^{0})^{0}=(\mathcal{L},\rightharpoonup ^{\ast })$.
\end{lem}

\begin{Beweis}
It is clear that $\rightharpoonup ^{\ast }$ is an $S-$action on $\mathcal{L}$%
. For all $s\in S$ and all $x\in L$ we have
\begin{equation}
s({\rightharpoonup ^{0}})^{0}\text{ }x=(s\rightharpoonup ^{0}1^{0})\vee
^{0}x=((s\rightharpoonup 1)\vee 0)\wedge x=(s\rightharpoonup 1)\wedge
x=s\rightharpoonup ^{\ast }x.
\end{equation}$\blacksquare$
\end{Beweis}

\begin{rems}\label{rems-new}
Let $(\mathcal{L},\rightharpoonup )=(L,\wedge,\vee ,0,1)$ a bounded lattice with an $S$-action.

\begin{enumerate}
\item $0$ is prime if and only if $1$ is first.

\item $SH(\mathcal{L})\subseteq Spec^p(\mathcal{L}^0)$.

\item If $(\mathcal{L},\rightharpoonup )$ is multiplication, then $$Spec^{psi}(\mathcal{L}) = SH(\mathcal{L})=Spec^{p}(\mathcal{L}^{0})$$.

Assume that $(\mathcal{L},\rightharpoonup )$ is multiplication. The first equality follow
from the definitions.

Let $x\in Spec^{p}(\mathcal{L}^{0})$. Suppose that $x\leq y\vee z$ for some $y,z\in L.$
Since $(\mathcal{L},\rightharpoonup )$ is multiplication, $%
y=s\rightharpoonup 1$ for some $s\in S$, and so $x\leq (s\rightharpoonup
1)\vee z$, i.e. $s\rightharpoonup ^{0}z\leq ^{0}x$. Since $x\in Spec^{p}(%
\mathcal{L}^{0}),$ we have $s\rightharpoonup ^{0}1^{0}\leq ^{0}x$ or $z\leq
^{0}x$ and so $x\leq s\rightharpoonup 1 = y$ or $x\leq z$. So, $Spec^{p}(%
\mathcal{L}^{0})\subseteq SH(\mathcal{L}).$ The inverse inclusion follows by
(2).

\item \index{pseudo strongly irreducible} $x\in L\backslash \{1\}$ is prime in $(\mathcal{L},\rightharpoonup
^{\ast })$ if and only if $x$ is pseudo strongly irreducible in $(\mathcal{L},\rightharpoonup )$.

\item $x\in L\backslash \{1\}$ is coprime in $(\mathcal{L},\rightharpoonup )$
if and only if $x$ is coprime in $(\mathcal{L},\rightharpoonup ^{\ast })$.

\item $x\in L\backslash \{0\}$ is first if and only if $0$ is prime in $%
[0,x].$

$(\Rightarrow )\ $Let $x\in L\backslash \{0\}$ be first. Observe that the
maximum element in the sublattice $[0,x]$ is $x$. Suppose that $%
s\rightharpoonup y=0$ for some $y\leq x$. Since $x$ is first, $y=0$ or $%
s\rightharpoonup x=0.$ So $0$ is prime in $[0,x].$

$(\Leftarrow )$ Let $0$ be prime in $[0,x].$ Suppose that $s\rightharpoonup
y=0$ for some $y\leq x.$ Since $y\in \lbrack 0,x]$, we have $y=0$ or $s\rightharpoonup
x=0$ as $x$ is the maximum element of $[0,x].$

\item $x\in L\backslash \{0\}$ is second if and only if $0$ is coprime in
the interval $[0,x].$
\end{enumerate}
\end{rems}

The notion of top-lattices was introduced by Abuhlail and Lomp \cite{AbuC}:

\begin{punto}
Let $(\mathcal{L},\rightharpoonup )=(L,\wedge,\vee ,0,1)$ a complete lattice and $X\subseteq L\backslash \{1\}.$
For $a\in L,$ we define the \textit{variety} of $a$\index{variety} as $V(a):=\{p\in X\mid a\leq p\}$ and set $V(\mathcal{L}%
):=\{V(a)\mid a\in L\}.$ Indeed, $V(\mathcal{L})$ is closed under arbitrary
intersections (in fact, $\bigcap_{a\in A}V(a)=V(\bigvee_{a\in A}(a))$ for
any $A\subseteq L$). The lattice $\mathcal{L}$ is called $X$\emph{-top}
(or a \textit{topological lattice} iff $V(\mathcal{L})$ is closed under finite
unions.
\end{punto}

Many results in the literature for prime, coprime, second, first, and other
types of spectra of submodules of a module can be generalized to a top-lattices
with actions from posets. For example, we have the following generalization of \cite[Theorem 3.5]{MMS1997}.

\begin{lemma}
Let $(\mathcal{L},\rightharpoonup )$ be a complete lattice with an action from a poset $S$. If $\mathcal{L}$ is multiplication, then $\mathcal{L}$ is $Spec^{p}(\mathcal{L})$-top.
\end{lemma}
\begin{Beweis} This follows from the fact that we have $V(s\rightharpoonup 1)\cup V(y)=V((s\rightharpoonup 1)\wedge y)$ for all $s\in S$ and $y\in L$. Indeed, by definition of prime elements and the axioms of the $S$-action, and noting that $V(-)$ is an order reversing map, we have:
\begin{equation*}
V(s\rightharpoonup y)\subseteq V(s\rightharpoonup 1)\cup V(y)\subseteq
V((s\rightharpoonup 1)\wedge y)\subseteq V(s\rightharpoonup y)\blacksquare.
\end{equation*}
\end{Beweis}

\begin{defn}
\label{quot}Let $\mathcal{L}=(L,\wedge ,\vee )$ be a lattice. Let $x,y,z\in
L $, with $x\leq y$ and $x\leq z$. We define $y\sim z$ iff for all $%
y^{\prime }\leq y$, there exists $z^{\prime }\leq z$ such that $y^{\prime
}\vee x=z^{\prime }\vee x$, and for all $z^{\prime }\leq z$, there exists $%
y^{\prime }\leq y$ such that $y^{\prime }\vee x=z^{\prime }\vee x$. It is
clear that $\sim $ is an equivalence relation. Denote the equivalence class
of $y\geq x$ by $y/x$, and define
\begin{equation*}
L/x:=\{y/x\mid y\in L\text{ and }x\leq y\}.
\end{equation*}%
Define $y/x\leq ^{q}z/x$ iff for all $y^{\prime }\leq y$, there exists $%
z^{\prime }\leq z$ such that $y^{\prime }\vee x=z^{\prime }\vee x$.
Then $\mathcal{L}/x=(L/x,\wedge ^{q},\vee ^{q})$ is a lattice,
called the \emph{quotient lattice}, where the meet
$\wedge ^{q}$ and the join $\vee ^{q}$ on $L/x$ are defined by:
\begin{equation*}
y/x\wedge ^{q}z/x:=(y\wedge z)/x\text{ and }y/x\vee ^{q}z/x:=(y\vee z)/x.
\end{equation*}%
If $\mathcal{L}=(L,\wedge ,\vee ,0,1)$ is a complete lattice, then
$\mathcal{L}/x=(L/x,\wedge ^{q},\vee ^{q})$ is a complete lattice, where
\begin{equation}
\bigwedge_{i\in A}^{q}(x_{i}/x)=(\bigwedge_{i\in A}x_{i})/x\text{ and }%
\bigvee_{i\in A}^{q}(x_{i}/x)=(\bigvee_{i\in A}x_{i})/x).
\end{equation}%
\end{defn}

\begin{rem}
\label{quotient action}Let $(\mathcal{L},\rightharpoonup )$ a lattice with an $S$-action. Define for all $s\in S$
and $y/x\in \mathcal{L}/x$:
\begin{equation}
s\rightharpoonup ^{q}y/x=(s\rightharpoonup y)\vee x
\end{equation}%
Then $(\mathcal{L}/x,\rightharpoonup ^{q})$ is a lattice with an $S$-action.
\end{rem}

\begin{thm}
\index{coprime element}
\index{second element}
\index{$S-$action}
\label{Proposition 5.8}Let $(\mathcal{L},\rightharpoonup )=(L,\wedge ,\vee ,0,1)$ a complete lattice with an $S$%
-action.

\begin{enumerate}
\item $Spec^c(\mathcal{L}) = Spec^s(\mathcal{L}^0)$.

\item $Spec^c(\mathcal{L}^0) = Spec^s(\mathcal{L}^*)$.

\item If $x\in L\backslash \{1\}$ is prime, then
\begin{equation*}
Spec^{f}(\mathcal{L}/x)=(\mathcal{L}/x)\backslash \{x/x\}.
\end{equation*}

\item Assume that the following additional condition is satisfied for
our action:
\begin{eqnarray}
s &\rightharpoonup &(y\vee z)=s\rightharpoonup y\vee s\rightharpoonup z%
\text{ for all }s\in S\text{ and }y,z\in L  \label{C1}
\end{eqnarray}

Then $x\in L\backslash \{1\}$ is prime $\Leftrightarrow Spec^{f}(\mathcal{L}%
/x)=(\mathcal{L}/x)\backslash \{x/x\}$.
\end{enumerate}
\end{thm}

\begin{Beweis}
\begin{enumerate}
\item $p\in Spec^{c}(\mathcal{L})\Leftrightarrow s\rightharpoonup 1\leq p$
or $(s\rightharpoonup 1)\vee p=1$ for all $s\in S$\newline
\hphantom{ $p\in Spec^c(\mathcal{L})$} $\Leftrightarrow s\rightharpoonup
1\vee p=p$ or $s\rightharpoonup ^{0}p=0^{0}$  for all $s\in S$ \newline
\hphantom{ $p\in Spec^c(\mathcal{L})$} $\Leftrightarrow s\rightharpoonup
^{0}p=p$ or $s\rightharpoonup ^{0}p=0^{0}$  for all $s\in S$ \newline
\hphantom{ $p\in Spec^c(\mathcal{L})$} $\Leftrightarrow p\in Spec^{s}(%
\mathcal{L}^{0})$.

\item $p\in Spec^{c}(\mathcal{L}^{0})\Leftrightarrow s\rightharpoonup
^{0}1^{0}\leq ^{0}p$ or $(s\rightharpoonup ^{0}1^{0})\vee ^{0}p=1^{0}$.%
\newline
\hphantom{$p\in Spec^c(\mathcal{L}^0)$} $\Leftrightarrow (s\rightharpoonup
1)\vee 0\geq p$ or $((s\rightharpoonup 1)\vee 0)\wedge p=0$   for all $s\in S$ \newline
\hphantom{$p\in Spec^c(\mathcal{L}^0)$} $\Leftrightarrow (s\rightharpoonup
1)\wedge p=p$ or $(s\rightharpoonup 1)\wedge p=0$   for all $s\in S$ \newline
\hphantom{$p\in Spec^c(\mathcal{L}^0)$} $\Leftrightarrow s\rightharpoonup
^{\ast }p=p$ or $s\rightharpoonup ^{\ast }p=0$   for all $s\in S$ \newline
\hphantom{$p\in Spec^c(\mathcal{L}^0)$} $\Leftrightarrow p\in Spec^{s}(%
\mathcal{L}^{\ast })$.

\item Let $x\in L\backslash \{1\}$ be prime. \textbf{Claim:} $y/x\in
\mathcal{L}/x$ is first.

Let $s\rightharpoonup ^{q}z/x=x/x$ and $z/x\leq ^{q}y/x$ and suppose that $%
z/x\nleq x/x$. Then $((s\rightharpoonup z)\vee x)/x=x/x$. It follows that $%
((s\rightharpoonup z)\vee x)=x$, and hence $((s\rightharpoonup z)\leq x$.
Since $x$ is prime, $((s\rightharpoonup 1)\leq x$ or $z\leq x$. But $z\leq x$
implies that $z=x$, and so $z/x=x/x$. Therefore, $((s\rightharpoonup 1)\leq
x $, and so $(s\rightharpoonup 1)\vee x=x$. Hence $s\rightharpoonup
^{q}1/x=x/x $. Therefore, $s\rightharpoonup ^{q}y/x=x/x$.

\item Assume that the additional condition (\ref{C1}) is satisfied and that $Spec^{f}(\mathcal{L}/x)=(\mathcal{L}/x)\backslash
\{x/x\} $. \textbf{Claim:} $x$ is prime in $\mathcal{L}$.

Suppose that $s\rightharpoonup y\leq x$ and $y\nleq x$. Since $%
s\rightharpoonup y\leq x,$ we have $(s\rightharpoonup y)\vee x=x$. It
follows by (\ref{C1}) that $s\rightharpoonup (y\vee x)=s\rightharpoonup
y\vee s\rightharpoonup x$. Since $s\rightharpoonup x\leq x$, we have
\begin{equation*}
s\rightharpoonup (y\vee x)=s\rightharpoonup y\vee s\rightharpoonup x\leq
(s\rightharpoonup y)\vee x=x.
\end{equation*}
Therefore $(s\rightharpoonup (y\vee x)\vee x)/x=x/x$, whence $%
s\rightharpoonup ^{q}(y\vee x)/x=x/x$. But $1/x$ is first in $\mathcal{L}/x,$
whence $(y\vee x)/x=x/x$ or $s\rightharpoonup ^{q}1/x=x/x$. Notice that $%
(y\vee x)/x=x/x$ cannot happen as $y\nleq x$. Thus $s\rightharpoonup
^{q}1/x=x/x$. Whence $s\rightharpoonup 1\vee x=x,$ i.e. $s\rightharpoonup
1\leq x$.$\blacksquare$
\end{enumerate}
\end{Beweis}

\begin{rem}
Let $(\mathcal{L},\rightharpoonup )=(L,\wedge,\vee ,0,1)$ a complete lattice with an $S$-action. Since $Spec^{c}(\mathcal{L%
})=Spec^{s}(\mathcal{L}^{0})$ by \ref{Proposition 5.8} (2), the result on
the second spectrum can be dualized to the coprime spectrum.
\end{rem}


\section{PS-Hollow Representation}

\hspace{1cm} Throughout this Section, $R$ is a commutative ring with unity and $M$ is a non-zero $R$-module. We consider the poset $\mathcal{I}=(Ideal(R),\subseteq)$ of ideals of $R$ acting on the lattice $\mathcal{L} = Sub_R(M)$ of $R$-submodules of $M$ in the canonical way. We say that a proper $R$-submodule of $M$ is irreducible (resp. strongly irreducible, pseudo strongly irreducible, prime, coprime) iff it is so as an element of $Sub_R(M)$. On the other hand, we say that a non-zero $R$-submodule of $M$ is hollow (resp. strongly hollow, pseudo strongly hollow, second, first) iff it is so as an element of $Sub_R(M)$. For such notions for modules one might consult \cite{Abu2011-CA}, \cite{Abu2011-TA}, \cite{Abu2015}, \cite{Yas1995}, \cite{Yas2001},\cite{Wij2006}).

In \cite{AH-sr}, we introduced and investigated modules attaining \emph{second representations}, i.e. modules which are finite sums of \emph{second submodules} (see \cite{Ann2002}, \cite{Bai2009}). Since every second submodule is secondary, modules with secondary representations can be considered as generalizations of such modules. Secondary modules can be considered, in some sense, as dual to those of primary submodules.

In this section, we consider modules with \emph{pseudo strongly hollow representations}, i.e. which are finite sums of \emph{pseudo strongly hollow submodules}. Assuming suitable conditions, we prove existence and uniqueness theorems for modules with such representations (called \emph{PS-hollow representable modules}). This work is inspired by the theory of primary and secondary decompositions of modules over commutative rings (e.g. Ann2002).

\begin{punto}
\index{primary submodule}
\index{primary decomposition}
\index{minimal primary submodule} A proper $R$-submodule $N\lneq M$ is
called \emph{primary} \cite{AK2012} iff whenever $rx\in N$ for some $r \in R$ and $x\in M$, either $x\in N$
or $r^{n}M\subseteq N$ for some $n\in \mathbb{N}$. We say that $M_R$ has a \emph{primary decomposition} \cite%
{AK2012} iff there are primary submodules $N_{1},\cdots ,N_{n}$ of $M$ with $M=\bigcap_{i=1}^{n}N_{i}.$

Dually, an $R$-submodule $N\leq M$ is said to be a \emph{secondary submodule} (\cite{K1973}, \cite{M1973}) iff for any $r\in R$
we have $rN=N$ or $r^{n}N=0$ for some $n\in \mathbb{N}$. An $R$-module $M$ has a {secondary representation} iff $M=\sum\limits_{i=1}^{n}N_{i},$ where $N_{1},\cdots ,N_{n}$ are secondary $R$-submodules of $M$.
\end{punto}

The notion of a primary submodule can be dualized in different ways. Instead of considering such notions, we consider the \emph{exact dual} of a pseudo strongly irreducible submodule (defined in \ref{psi}). Recall that, the pseudo strongly irreducible elements in $(Sub_R(M),\rightharpoonup)$ are exactly the prime elements in $(Sub_R(M),\rightharpoonup^*)$ (defined in \ref{psh}).

Strongly irreducible submodules (ideal) have been considered by several authors (e.g. \cite{HRR2002}, \cite{Ata2005}, \cite{Azi2008}). The dual notion of a \emph{strongly hollow submodule} was investigated by Abuhlail and Lomp in \cite{AbuC0}. In this section we consider the more general notion of a \emph{pseudo strongly hollow submodule}. For the convenience of the reader, we restate the definition  in the special case of the lattice $Sub_R(M)$.

\begin{defn}
We say that an $R$-submodule $N\leq M$ is \emph{pseudo strongly hollow submodule }%
\index{pseudo strongly hollow submodule} (or \emph{PS-hollow}%
\index{PS-hollow submodule} for short) iff for any ideal $I\leq R$ and any $R
$-submodule $L\leq M,$ we have%
\begin{equation}
N\subseteq IM+L\Rightarrow N\subseteq IM%
\text{ or }N\subseteq L.  \label{PSH}
\end{equation}%
We say that $_{R}M$ is a \emph{pseudo strongly hollow module} (or PS-hollow\emph{\ }%
for short) iff $M$ is a PS-hollow submodule of itself, that is, $M$ is
PS-hollow iff for any ideal $I\leq R$ and any $R$-submodule $L\leq M,$ we
have%
\begin{equation}
M=IM+L\Rightarrow M=IM\text{ or }M=L.  \label{PH}
\end{equation}
\end{defn}

\begin{ex}
Let $_{R}M$ be second. Every $R$-submodule $N\leq M$ is a PS-hollow
submodule of $M.$ Indeed, suppose that $N\subseteq IM+L$ for some $L\leq M$
and $I\leq R$. Since $_{R}M$ is second, either $IM=0$ whence $N\subseteq L,$
or $IM=M$ whence $N\subseteq IM$. In particular, every second module is a
PS-hollow module.
\end{ex}

\begin{rem}
It is clear that any strongly hollow submodule of $M$ is PS-hollow; the
converse holds in case $_{R}M$ is multiplication.
\end{rem}

\begin{ex}
\item \label{Example 5.2}

\begin{enumerate}
\item There exists an $R$-module $M$ which is not multiplication but all of
its PS-hollow submodules are strongly hollow. Consider the Pr\"{u}fer group $M=\mathbb{Z}({p^{\infty }})$ as a $\mathbb{Z}$-module. Notice that $_{\mathbb{Z}}M$ is not a multiplication module, however every $\mathbb{Z}$-submodule of $M$ is strongly hollow).

\item A PS-hollow submodule $N\leq M$ need not be hollow. Consider $M=%
\mathbb{Z}_{2}[x]$ as a $\mathbb{Z}$-module. Set $N:=x\mathbb{Z}_{2}[x]$,
and $L:=(x+1)\mathbb{Z}_{2}[x]$. Then $N,L\lneq M$ and $M=L+N$ is PS-hollow
which is not hollow. Indeed, $x^{i}=x^{i-1}(x+1)-x^{i-2}(x)$ for all $i\geq
2 $ and $1=(x+1)-x$. On the other hand, $IM=M$ or $IM=0$ for
every $I\leq \mathbb{Z}$.
\end{enumerate}
\end{ex}

\begin{lem}
\label{Lemma 5.3}Let $N\leq M$ be a PS-hollow submodule. If $I$ is minimal
in $A:=\{I\leq R\mid N\subseteq IM\},$ then $I$ is a hollow
ideal of $R.$
\end{lem}

\begin{Beweis}
Let $I=J+K$ for some ideals $J,K\leq R.$ Notice that $N\subseteq
IM=(J+K)M=JM+KM$, whence $N\subseteq JM$ or $N\subseteq KM$, i.e. $J\in A$
or $K\in A$. By the minimality of $I,$ it follows that $J=I$ or $K=I$.
Therefore, $I$ is hollow.$\blacksquare$
\end{Beweis}

\begin{punto}
Let $N\leq M$ be a PS-hollow submodule and set%
\begin{equation*}
A_{N}:=\{I\leq R\mid N\subseteq IM\},\text{ }H_{N}:=Min(A)\text{ and }%
In(N):=\bigcap\limits_{I\in H_N}IM.
\end{equation*}%
Notice that $A_{N}$ is non-empty as $R\in A,$ while $H_{N}$ might be empty
and in this case $In(N):=M$ (however $H_{N}\neq \emptyset $ if $R$ is
Artinian). When $N$ is clear from the context, we drop it from the index of the above notations. We
say that $N$ is an $H$\emph{-PS-hollow submodule}%
\index{$H$-PS-hollow submodule} of $M.$ Every element in $H$ is called an
\emph{associated hollow ideal of }$M$%
\index{associated hollow ideal}. We write $Ass^{h}(M)$ to denote the set of
all associated hollow ideals of $M.$
\end{punto}

\begin{prop}
\index{PS-hollow submodule}
\label{Proposition 5.5}Let $R$ be an Artinian ring, $N$ and $L$ be
incomparable PS-hollow submodules of $M$ and $H\subseteq Ass^{h}(M)$. Then $%
N+L$ is $H$-PS-hollow if and only if $N$ and $L$ are $H$-PS-hollow.
\end{prop}

\begin{Beweis}
$(\Leftarrow )$ Let $N\leq M$ and $L\leq M$ be $H$-PS-hollow submodules.

\textbf{Claim 1:} $H_{N+L}=H_{N}=H.$

Consider $I\in H_{N}=H_{L}.$ Clearly, $I\in A_{N+L}.$ If $I\notin
H_{N+L}:=Min(A_{N+L})$, then there is $I^{\prime }\subsetneq I$ such that $%
N\subseteq N+L\subseteq I^{\prime }M$ which contradicts the minimality of $I$
in $A_{N}$.

Conversely, let $I\in H_{N+L}$. Clearly, $I\in A_{N}\cap A_{L}.$ If $I\notin
H_{N},$ then there is $I^{\prime }\in H_{N}=H_{L}$ with $I^{\prime
}\subseteq I$ and therefore $N+L\subseteq I^{\prime }M,$ whence $I=I^{\prime
}$ since $I^{\prime }\in A_{N+L}.$ Therefore, $H_{N+L}=H_{N}=H$.

\textbf{Claim 2:} $N+L$ is PS-hollow.

Suppose that $N+L\subseteq JM+K$ for some ideal $J\leq R$ and some submodule $%
K\leq M$. Then $N\subseteq N+L\subseteq JM+K$ and so $N\subseteq JM$ or $%
N\subseteq K$. Similarly $L\subseteq N+L\subseteq JM+K$ and so $L\subseteq
JM $ or $L\subseteq K$. Suppose that $N\subseteq JM$, whence there is $I\in
H $ such that $N\subseteq IM$ and $I\subseteq J$ (as $R$ is Artinian) and so
$L\subseteq IM\subseteq JM$. Therefore, either $N+L\subseteq JM$ or $%
N+L\subseteq K$. Hence $N+L$ is PS-hollow.

$(\Rightarrow )$ Assume that $N+L$ is $H$-PS-hollow. It is clear that $%
H_{N+L}\subseteq H_{L}$. Assume that $L\subseteq IM$. Then $N+L\subseteq
IM+L $ and $N+L\nsubseteq L$ as $N$ and $L$ are incomparable, whence $%
N+L\subseteq IM$ and so $H_{L}\subseteq H_{N+L}$. Therefore, $H_{L}=H_{N+L}$%
. Similarly, $H_{N}=H_{N+L}$.$\blacksquare$
\end{Beweis}

\begin{punto}
\index{DPS-hollow representable module}%
\index{SPS-hollow representable module}%
\index{minimal PS-hollow representation} We say that a module $M$ is \emph{%
PS-hollow representable}%
\index{PS-hollow representable module} iff $M$ can be written as a \emph{%
finite} sum of PS-hollow submodules. A module $M$ is called \emph{directly
PS-hollow representable}%
\index{directly PS-hollow representable} (or DPS-hollow
representable, for short) iff $M$ is a \emph{finite direct} sum of
PS-hollow submodules. A module $M$ is called \emph{semi-pseudo strongly
hollow representable} (or \emph{SPS-hollow representable}, for short) iff $M$
is a sum of PS-hollow submodules. We call $M=\sum\limits_{i=1}^{n}N_{i},$
where each $N_{i}$ is $H_{i}$-PS-hollow, a \emph{minimal PS-hollow
representation}%
\index{minimal PS-hollow representation} for $M$ iff the following
conditions are satisfied:

\begin{enumerate}
\item $In(N_{1}), In(N_{2}), .... , In(N_{n})$ are incomparable.

\item $N_{j}\nsubseteq \sum\limits_{i=1,i\neq j}^{n}N_{i}$ for all $j\in
\{1,\cdots ,n\}$.
\end{enumerate}

If such a minimal PS-hollow representation for $M$ exists, then we call each
$N_{i}$ a \emph{main PS-hollow submodule}
\index{main PS-hollow submodule}\emph{of} $M$ and the elements of $%
H_{1},H_{2},\cdots ,H_{n}$ are called \emph{main associated hollow ideals}
\index{main associated hollow ideals} of $M;$ the set of the main associated
hollow ideals of $M$ is dented by $ass^{h}(M)$.
\end{punto}

\begin{prop}
\index{minimal PS-hollow representation}
\label{Remark 5.6}(Existence of minimal PS-hollow representation) Let $R$ be an Artinian ring and suppose that  $In(N)$ is PS-hollow whenever $N$ is PS-hollow. Then every PS-hollow representable $R$-module has a minimal PS-hollow representation.
\end{prop}

\begin{Beweis}
Let $M=\sum\limits_{i\in A}K_{i}$, where $A$ is finite and $K_{i}$ is an $%
H_{i}$-PS-hollow submodule $\forall i\in A$.

\textbf{Step 1:} Remove the redundant submodules $K_{j}\subseteq
\sum\limits_{i\neq j}K_{i}.$ This is possible by the finiteness of $A$.

\textbf{Step 2:} Gather all submodules $K_{m}$ that share the same $H$ to
construct an $H$-PS-hollow $N\leq M$ as a sum of such $H$-PS-hollow
submodules (this is possible by Proposition \ref{Proposition 5.5}).

\textbf{Step 3:} If $In(K_{i})$ and $In(K_{j})$ are comparable for some $i,j\in \{1, 2, .., n\}$; say $In(K_{i}) \subseteq In(K_{j})$ then replace $K_i$ and $K_j$ in the representation by $In(K_{j})$.$\blacksquare$
\end{Beweis}

\begin{ex} Any vector space $V$ has a trivial minimal PS-hollow representation as it is a PS-hollow submodule of itself. Notice that $V$ is not necessarily multiplication.
\end{ex}

We provide an example of a module with a minimal PS-hollow representation that is \it{not} a strongly hollow representation.

\begin{ex} Consider $R:={\Bbb{Z}}_{pq}$ where $p$ and $q$ are distinct prime numbers and $M= {\Bbb{Z}}_{pq}[x]$. Notice that $M= pM+qM$ is a minimal PS-hollow representation while neither $pM$ nor $qM$ is strongly hollow. To see this, notice that $M=xM + {\Bbb{Z}}_{pq}$ but neither $pM\subseteq xM$ nor $pM\subseteq {\Bbb{Z}}_{pq}$ (similarly, neither $qM\subseteq xM$ nor $qM\subseteq {\Bbb{Z}}_{pq}$). Observe that $M$ is neither multiplication nor a vector space.
\end{ex}

\begin{rem}
\index{$H$-PS-hollow submodule}
\label{In(N) is H-PS-hollow}Let $R$ be Artinian and $N\leq M$ be an $H$%
-PS-hollow submodule. If $In(N)$ is PS-hollow, then $In(N)$ is $H$%
-PS-hollow. To show this, observe that for any ideal $I\leq R$, we have $%
N\subseteq IM$ if and only if there exists $I^{\prime }\in H$ such that $%
N\subseteq I^{\prime }M$ with $I^{\prime }\subseteq I$ (as $R$ is Artinian),
whence $In(N)\subseteq IM$ if and only if $N\subseteq IM$.
\end{rem}

\begin{lem}
\label{Lemma 5.7}Let $R$ be Artinian, $N\leq M$ be an $H$-PS-hollow submodule and $%
In(N)\leq L$ whenever $N\leq L\leq M$. Then $In(N)$ is $H$-PS-hollow.
\end{lem}

\begin{Beweis}
Let $K=In(N):=\bigcap\limits_{I\in H}IM$. Suppose that $K\subseteq JM+L$ for
some $J\leq R$ and $L\leq M.$ If $K\nsubseteq JM$, then $N\nsubseteq JM$ and
so $N\subseteq L$, whence $K\subseteq L$. Therefore $K$ is PS-hollow. Thus,
by the Remark \ref{In(N) is H-PS-hollow}, $In(N)$ is $H$-PS-hollow.$\blacksquare$
\end{Beweis}

\begin{ex}
If $R$ is Artinian, then every multiplication $R$-module $M$ satisfies the
conditions of Lemma \ref{Lemma 5.7} and so $In(N)$ is $H$-PS-hollow for
every $H$-PS-hollow $N\leq M$ (in fact, $In(N)=N$ in this case).
\end{ex}

\begin{rem}
\label{Example 5.8}Let $R$ be Artinian and $M$ a multiplication $R$-module.
It is easy to see that there is a unique minimal PS-hollow representation of
$M$ up to the order, \emph{i.e.} if $\sum\limits_{i=1}^{n}N_{i}=M=\sum%
\limits_{j=1}^{m}K_{j}$ are two minimal PS-representations such that each $%
N_{i}$ is $H_{i}$-PS-hollow and each $K_{j}$ is $H_{j}^{\prime }$-PS-hollow,
then $n=m$ and $\{N_{1},\cdots ,N_{n}\}=\{K_{1},\cdots ,K_{n}\}.$
\end{rem}

\begin{thm}

\label{Theorem 5.9}\index{minimal PS-hollow representation}(First uniqueness theorem of PS-hollow representation)
Let $R$ be Artinian and $\sum\limits_{i=1}^{n}N_{i}=M=\sum%
\limits_{j=1}^{m}K_{j}$ be two minimal PS-representations for $_{R}M$ such
that $N_{i}$ is $H_{i}$-PS-hollow for each $i\in \{1,\cdots ,n\}$ and $K_{j}$
is $H_{j}^{\prime }$-PS-hollow for each $j\in \{1,\cdots ,m\}$. Then $n=m,$ $%
\{H_{1},\cdots ,H_{n}\}=\{H_{1}^{\prime },\cdots,H_{n}^{\prime }\}$ and $In(N_{i})=In(K_{j})$
whenever $H_{i}=H_{j}^{\prime} $.
\end{thm}

\begin{Beweis}
Set $N_{i}^{\prime }=In(N_{i})$ and $K_{j}^{\prime }=In(K_{j})$ for $i\in
\{1,\cdots ,n\}$ and $j\in \{1,\cdots ,m\}$.

\textbf{Claim: }For any $i\in \{1,\cdots ,n\}$, there is $j\in
\{1,\cdots ,m\}$ such that $N_{i}^{\prime }=K_{j}^{\prime }$.

\textbf{Step 1:} Suppose that there exists some $i\in \{1,\cdots ,n\}$ for
which $N_{i}\nsubseteq K_{j}^{\prime }$ for all $j\in \{1,\cdots ,m\}$.
Then for any $j\in \{1,\cdots ,m\}$, there is $J_{j}^{\prime }\in
H_{j}^{\prime }$ such that $N_{i}\nsubseteq J_{j}^{\prime }M.$ But $%
N_{i}\subseteq M=\sum\limits_{j=1}^{n}K_{j}\subseteq
\sum\limits_{j=1}^{m}J_{j}^{\prime }M$, whence $N_{i}\subseteq J_{j}^{\prime
}M$ for some $j$ (a contradiction). So, $N_{i}\subseteq K_{j}^{\prime }$ for
some $j\in \{1,\cdots ,m\}.$

\textbf{Step 2:}\ We show that $N_{i}^{\prime }\subseteq K_{j}^{\prime }$.

Since $N_{i}\subseteq K_{j}^{\prime }$, we have $N_{i}\subseteq IM$ for all $%
I\in H_{j}^{\prime }$. Since $R$ is Artinian, there is a minimal ideal $%
J_{I}\leq I$ such that $N_{i}\subseteq J_{I}M$ and so
\begin{equation*}
N_{i}^{\prime }=In(N_{i})=\bigcap_{I\in H_{i}}IM\subseteq \bigcap_{I\in
H_{j}^{\prime }}J_{I}M\subseteq K_{j}^{\prime }.
\end{equation*}%
Similarly, for any $j\in \{1,\cdots ,m\}$, there is some $i\in
\{1,\cdots ,n\}$ such that $K_{j}^{\prime }\subseteq N_{i}^{\prime }$.
Therefore, for any $i\in \{1,\cdots ,n\}$, there is some $j\in
\{1,\cdots ,m\}$ such that $N_{i}^{\prime }=K_{j}^{\prime }$ as $%
N_{1}^{\prime },N_{2}^{\prime },\cdots ,N_{n}^{\prime }$ are incomparable.

\textbf{Claim:} $H_{i}=H_{j}^{\prime }$ whenever $N_{i}^{\prime
}=K_{j}^{\prime }.$

Let $N_{i}^{\prime }=K_{j}^{\prime }.$ Pick any $I\in H_{i}$. Then $%
N_{i}\subseteq IM$, whence $K_{j}^{\prime }=N_{i}^{\prime }\subseteq IM$.
Since $R$ is Artinian, there is a minimal ideal $I^{\prime }\in
H_{j}^{\prime }$ such that $I^{\prime }\leq I$, and therefore $I^{\prime }=I$
as $I$ is minimal with respect to $N_{i}\subseteq IM$. Hence $H_{i}\subseteq
H_{j}^{\prime }$. One can prove similarly that $H_{j}^{\prime }\subseteq
H_{i}.$ So, $H_{i}=H_{j}^{\prime }$.$\blacksquare$
\end{Beweis}

\begin{thm}
\label{Theorem 5.10}\index{minimal PS-hollow representation}(Second uniqueness theorem of PS-hollow representation)
Let $R$ be Artinian, $M$ be an $R$-module with two minimal PS-hollow
representations $\sum\limits_{i=1}^{n}N_{i}=M=\sum\limits_{j=1}^{n}K_{j}$
with $N_{i}$ is $H_{i}$-PS-hollow for each $i\in \{1,\cdots ,n\}$ and $K_{j}$
is $H_{j}$-PS-hollow for each $j\in \{1,\cdots ,n\}$. If $H_{m}$ is minimal
in $\{H_{1},H_{2},\cdots ,H_{n}\}$, then either $N_{m}=K_{m}$ or $In(N_{m})$
is not PS-hollow.
\end{thm}

\begin{Beweis}
Let $H_{m}$ be minimal in $\{H_{1},H_{2},\cdots ,H_{n}\}$ such that $%
In(N_{m})$ is PS-hollow. For any $j\neq m$, there is $I_{j}\in
H_{j}\backslash H_{m}.$ But $\sum\limits_{j\neq m}I_{j}M+N_{m}=M$ and so $%
In(N_{m})\subseteq \sum\limits_{j\neq m}I_{j}M+N_{m}$. Since $I_{j}\in
H_{j}\backslash H_{m}$, it follows that $In(N_{m})\nsubseteq I_{j}M$ for all
$j\in \{1,\cdots ,n\}\backslash \{m\}$ and so $In(N_{m})\subseteq N_{m}$,
whence $In(N_{m})=N_{m}$. One can prove similarly that $In(K_{m})=K_{m}.$ It
follows that%
\begin{equation*}
N_{m}=In(N_{m})\overset{%
\text{Theorem }\ref{Theorem 5.9}}{=}In(K_{m})=K_{m}.
\end{equation*}$\blacksquare$
\end{Beweis}

\begin{cor}
\label{Corollary 5.11}Let $R$ be Artinian and $\sum\limits_{i=1}^{n}N_{i}=M=%
\sum\limits_{i=1}^{n}K_{i}$ be two minimal PS-hollow representations of $%
_{R}M$ such that $N_{i}$ is $H_{i}$-PS-hollow for $i\in \{1,\cdots ,n\}$
and $K_{i}$ is $H_{i}$-PS-hollow for $i\in \{1,\cdots ,n\}$. If $In(N)$ is
PS-hollow whenever $N$ is a main PS-hollow submodule of $M,$ then $%
N_{i}=K_{i}$ for all $i\in \{1,\cdots ,n\}$.
\end{cor}

\begin{Beweis}
Apply Theorem \ref{Theorem 5.10} and observe that $H_{i}$ is minimal in $%
\{H_{1},H_{2},\cdots ,H_{n}\}$ for each $i\in \{1,\cdots ,n\}$ as $%
In(N_{i})$ is PS-hollow: otherwise, $H_{j}\subsetneq H_{i}$ for some $i\neq
j $ and $In(N_{j})$ can replace $N_{i}+N_{j}$ whence $\sum%
\limits_{i=1}^{n}N_{i}$ is not minimal (a contradiction).$\blacksquare$
\end{Beweis}

\begin{punto}
We say that an $R$-module $M$ is\emph{\ pseudo distributive}%
\index{pseudo distributive module} iff for all $L,N\leq M$ and every $I\leq
R $ we have
\begin{equation}
L\cap (IM+N)=(L\cap IM)+(L\cap N).
\end{equation}%
Every distributive $R$-module is indeed pseudo distributive. The two notions
coincide for multiplication modules.
\end{punto}

\begin{ex}
A pseudo distributive module need not be distributive. Consider $M:=\mathbb{Z%
}_{2}[x]$ as a $\mathbb{Z}$-module. Let $N:=xM$, $L:=(x+1)M$ and $K=\mathbb{Z%
}_{2}$. Then $N,L,K\leq M$ are $R$-submodules and
\begin{equation*}
(K\cap L)+(K\cap N)=0\neq K=K\cap (L+N).
\end{equation*}%
Notice that $M$ is pseudo distributive as $IM=0$ or $IM=M$ for every $I\leq
R $.
\end{ex}

\begin{rem}
\index{directly hollow representable}
\index{hollow representable}
  Assume that $M$ is a (directly) hollow representable module for which every maximal hollow is PS-hollow. Then $M$ is (directly) PS-hollow representable.
\end{rem}

In \cite{AH-sr}, we introduced the notion of \emph{s-lifting modules}:

\begin{punto}
Recall that an $M$ is a \emph{lifting} $R$-module iff any $R$%
-submodule $N\leq M$ contains a direct summand $X\leq M$ such that $N/X$ is
small in $M/X$ (e.g. \cite[22.2]{JCNR}). We call $_{R}M$ \emph{s-lifting} iff $_{R}M$ is lifting
and every \emph{maximal} hollow submodule of $M$ is second.
\end{punto}

\begin{prop}
\label{hollow-PS}
\index{s-lifting}

\begin{enumerate}
\item If $_{R}M$ is pseudo distributive, then every hollow submodule of $M$
is PS-hollow.

\item If $_{R}M$ is s-lifting, then every maximal hollow submodule of $M$ is
PS-hollow.
\end{enumerate}
\end{prop}

\begin{Beweis}
\begin{enumerate}
\item Let $M$ is pseudo distributive. Let $N\leq M$ be hollow. Suppose that $%
N\subseteq IM+L$ , whence $N=(IM+L)\cap N=(IM\cap N)+(L\cap N)$ as $M$ is pseudo
distributive. Since $N$ is hollow, $N=IM\cap N$ or $N=L\cap N$, therefore $%
N\subseteq IM$ or $N\subseteq L$. So, $N$ is PS-hollow.

\item Let $_{R}M$ be s-lifting. Suppose that $K\leq M$ is a maximal hollow
submodule of $M$ and that $K\leq IM+L$. Since $M$ is s-lifting, there exists
$K^{\prime }\subseteq K$ and $N\leq M$ such that $K^{\prime }\oplus N=M$ and
$K/K^{\prime }$ is small in $M/K^{\prime }$.

\textbf{Case 1: }$K^{\prime }=0$: i.e. $M=N$. Since $K$ is second, we
have $K=IK\subseteq IN=IM$.

\textbf{Case 2: }$K^{\prime }\neq 0$: We claim that $K=K^{\prime }$. To
prove this, let $x\in K$. Then there are $y\in K^{\prime }$ and $z\in N$
such that $x=y+z$. But $y\in K$, whence $z\in K$. Therefore, $K\subseteq
K^{\prime }\oplus (K\cap N)$, but $K$ hollow implies that $K=K^{\prime }$ or
$K=K\cap N$. But $K^{\prime }\neq 0$, whence $K=K^{\prime }$; otherwise, $%
K^{\prime }\cap N\neq 0$. Therefore, $M=K\oplus N$. Now, it is easy to show
that
\begin{equation*}
IM+L\leq (IM\cap K+L\cap K)\oplus (IM\cap N+L\cap N),
\end{equation*}%
and so
\begin{equation*}
K\leq (IM\cap K+L\cap K)\oplus (IM\cap N+L\cap N),
\end{equation*}%
whence $K\leq IM\cap K+L\cap K$. Since $IM\cap K+L\cap K\leq K$, it follows
that $K=IM\cap K+L\cap K$ and so $K=IM\cap K$ or $K=L\cap K$ which implies
that $K\leq IM$ or $K\leq L$.$\blacksquare$
\end{enumerate}
\end{Beweis}

\begin{exs}
\label{Example 5.14}
\index{PS-hollow representable}
\index{directly PS-hollow representable}
\index{directly hollow representable}
\index{hollow representable}
\begin{enumerate}
\item Every (directly) hollow representable pseudo distributive module is (directly) PS-hollow representable.
\item Every s-lifting module with finite hollow  dimension is directly PS-hollow representable.
\item The $\mathbb{Z}$-module $M=\mathbb{Z}_{n}$ is PS-hollow representable.
To see this, consider the prime factorization $n=p_{1}^{m_{1}}\cdots
p_{k}^{m_{k}}$, and set $n_{i}=%
\frac{n}{p_{i}^{m_{i}}}$ for $i\in \{1,\cdots ,k\}.$ Then $%
M=\sum\limits_{i=1}^{k}(n_{i})$ is a minimal PS-hollow representation for $M$%
, and $(n_{i})$ is $H_{i}$-PS-hollow where $H_{i}=\{(n_{i})\}$ for $i\in
\{1,\cdots ,k\}$.

\item The $\mathbb{Z}$-module $M=\mathbb{Z}_{12}$ is PS-hollow representable
($M=4Z_{12}+3Z_{12}$), but $M$ is not second representable. Observe that $M$
is not semisimple and is even not s-lifting as $3Z_{12}\leq \mathbb{Z}_{12}$
is a maximal hollow $\mathbb{Z}$-subsemimodule but not second.

\item Any Noetherian semisimple $R$-module is directly PS-hollow
representable.

\item Any Artinian semisimple $R$-module is directly PS-hollow representable.
\end{enumerate}
\end{exs}

\begin{lem}
\label{Lemma 5.15}Let $N\leq M$ be an $H$-PS-hollow submodule such that
every non-small submodule $K$ of $M$ is of the form $JM$ for some ideal $%
J\leq R$. Every non-small submodule $K\leq N$ is $H$-PS-hollow submodule; Moreover, for any ideal $I\leq R$, we have: $K\subseteq IM$ if and only if $N\subseteq IM$.
\end{lem}

\begin{Beweis}
Let $N\leq M$ be an $H$-PS-hollow submodule and $K\leq N$ be a non-small submodule.
Suppose that $K\subseteq IM+L$ and $K\nsubseteq L$. Notice that $N\nsubseteq L$.
Since $K$ is not small in $N$, there is a proper submodule $K^{\prime }$ of $N$ such that
$N=K+K^{\prime } \subseteq IM+L+K^{\prime }$. If $N\subseteq L+K^{\prime }$, then $K^{\prime }=JM$ for some $J\leq R$ (notice that $K^{\prime }$ not small in $N$) and therefore $N\subseteq K^{\prime }$ (a contradiction).
Hence, $N\subseteq IM$ and so $K\subseteq IM$, whence $K$ is PS-hollow.

\textbf{Claim:} $A_H = A_K$.

Assume that $K\subseteq IM$ for some $I\leq R$. Then $%
N=K+K^{\prime }\subseteq IM+K^{\prime }$. Since $N$ is PS-hollow and $%
K^{\prime }\neq N$, we have $N\subseteq IM$.
\end{Beweis}

\begin{ex}
\label{Example 5.16}Consider $M=\mathbb{Z}_{12}$ as a $\mathbb{Z}$-module.
Then $K_{1}=3\mathbb{Z}_{12}$ and $K_{2}=4\mathbb{Z}_{12}$ satisfy the
assumptions of Lemma \ref{Lemma 5.15}. Notice that $_{\mathbb{Z}}M$ is not
semisimple.
\end{ex}

\begin{punto}
\index{comultiplication module}
  A module $_{R}M$ is called \emph{comultiplication} \cite{Abu2011-TA} iff
for every submodule $K\leq M$, we have $K=(0:_{M}(0:_{R}K)).$
\end{punto}

\begin{thm}
\index{PS-hollow submodule}
\index{second submodule}
\index{multiplication module}
\index{comultiplication module}
\label{Theorem 5.17}Let $_{R}M$ be semisimple, $B$ the set of all maximal
second submodules of $M$, and assume that $Ann(M)\neq \bigcap\limits_{K\in
B\backslash \{N\}}Ann(K)$ for any $N\in B$. The following conditions are
equivalent:

\begin{enumerate}

\item $_{R}M$ is multiplication.

\item Every PS-hollow submodule of $M$ is simple.

\item Every second submodule of $M$ is simple.

\item $_{R}M$ is comultiplication.
\end{enumerate}
\end{thm}

\begin{Beweis}
Let $M=\bigoplus\limits_{S\in A}S$, where $S$ is a simple submodule of $M$ for
all $S\in A$.

$(1)\Rightarrow (2)$: Assume that $_{R}M$ is multiplication. Suppose that there is an
$H$-PS-hollow submodule $N\leq M$, which is not simple. Then $N$ contains properly a
simple submodule $S^{\prime }\in A$. Since $S^{\prime }$ is not small in $N$,
Lemma \ref{Lemma 5.15} implies that $S^{\prime }$ is $H$-PS-hollow. But
there is another simple submodule $S^{\prime \prime }$ of $N$ (as $N$ is not
simple). Let $I:=Ann(S^{\prime \prime }).$ It follows that $S^{\prime}\subseteq IM$ while
$N\nsubseteq JM$ (which contradicts Lemma \ref{Lemma 5.15}).

\vspace*{0.5cm}

$(2)\Rightarrow (3)$: Assume that every PS-hollow submodule of $M$ is simple.

\textbf{Claim:} Every second submodule of $M$ is PS-hollow, whence simple.

Let $N = \bigoplus\limits_{i\in A}S_{i}$ be a second submodule of $M$ and
suppose that $N\subseteq IM+L$ for some ideal $I$ of $R$ and some $R$-submodule $N$ of $M$.

{\it Case 1}: $I\subseteq Ann(N)$. In this case, $N\cap IM=0$, and it follows that $N\subseteq L$.

{\it Case 2}: $I\nsubseteq Ann(N)$. Since $N$ is second, $N=IN\subseteq IM$.

\vspace*{0.5cm}

$(3)\Rightarrow (1)$: Assume that every second submodule of $M$ is simple.
Consider a submodule $K=\bigoplus\limits_{S\in C\subseteq A}S$ of $M$ and
set $I:=\bigcap\limits_{S\in A\backslash C}Ann(S)$. Notice that $K=IM$,
otherwise, $I\subseteq Ann(S)$ for some $S\in C$ whence $Ann(M)=\bigcap_{S%
\in A\backslash \{S\}}Ann(S)$ (a contradiction).
Since $K$ is an arbitrary submodule of $M$, we conclude that $M_R$ is multiplication.

\vspace*{0.5cm}

$(3)\Rightarrow (4)$: Assume that every second submodule of $M$ is simple.
Consider a submodule $K=\bigoplus\limits_{S\in C\subseteq A}S$ of $M$ and
set $I:=(0:_{R}K)$. Suppose that $(0:_{M}I)\neq K$, whence there is a simple
submodule $S^{\prime }\leq M$ with $S^{\prime }\cap K=0$ and $I\subseteq
Ann(S^{\prime })$ which is not allowed by our assumption as it would yield $%
Ann(M)=\bigcap\limits_{S\in B}Ann(S)=\bigcap\limits_{S\in B\backslash
\{S^{\prime }\}}Ann(S)$ (a contradiction to the assumption).

\vspace*{0.5cm}

$(4)\Rightarrow (3)$: Let $_{R}M$ be comultiplication. Let $K\leq M$ be second.
For any simple $S\leq K$ we have%
\begin{equation}
K=(0:_{M}(0:_{R}K))=(0:_{M}(0:_{R}S))=S,
\end{equation}%
i.e. $_{R}K$ is simple.$\blacksquare$
\end{Beweis}

\begin{ex}
\label{Example 5.18}Consider the $\mathbb{Z}$-module $M=\prod\limits_{i=1}^{%
\infty }\mathbb{Z}_{p_{i}p_{i}^{\prime }}$, where $p_{i}$
and $p_{i}^{\prime }$ are primes and $p_{i}\neq p_{j}$, $p_{i}^{\prime }\neq
p_{j}^{\prime }$ for all $i\neq j\in \mathbb{N}$ and $p_{i}^{\prime }\neq
p_{j}$ for any $i$ and $j$. Let the simple $\mathbb{Z}$-modules $K_{p_{i}}$
and $K_{p_{i}^{\prime }}$ be such that $(0:K_{p_{i}})=(p_{i})$ and $%
(0:K_{p_{i}^{\prime }})=(p_{i}^{\prime })$, so
\begin{equation*}
M=\bigoplus\limits_{i=1}^{\infty }K_{p_{i}}\oplus \bigoplus_{i=1}^{\infty
}K_{p_{i}^{\prime }}.
\end{equation*}%
Every second $\mathbb{Z}$-submodule of $M$ is simple, while $_{\mathbb{Z}}M$
is not multiplication. Notice that the assumption on $Ann(M)$ in Theorem \ref%
{Theorem 5.17} is not satisfied for this $\mathbb{Z}$-module, which shows
that this condition cannot be dropped.
\end{ex}

Recall from \cite{AH-sr} that an $R$-module $M$ is \emph{second representable} iff $M=\sum\limits_{i=1}^{n}K_{i}$, where $K_{i}$ is a second $R$-submodule of $M$ for all $i=1,\cdots ,n$. If this \emph{second representation} is minimal, the set of \emph{main second attached primes} of $M$ is given by $att^{s}(M)=\{Ann(K_{i})\mid i=1,\cdots ,n\}$.

\begin{cor}
\label{Corollary 5.19}If $_{R}M$ is semisimple second representable with $att^{s}(M)=Min(att^{s}(M))$. The following are equivalent:

\begin{enumerate}

\item $M$ is multiplication.

\item Every PS-hollow submodule of $M$ is simple.

\item Every second submodule of $M$ is simple.

\item $M$ is comultiplication.
\end{enumerate}
\end{cor}

\begin{Beweis}
Since $M$ is second representable, the set $B$ defined in Theorem \ref%
{Theorem 5.17} is finite. Since $Ann(S_{i})$ is prime for every $i\in A$ and
$att^{s}(M)=Min(att^{s}(M))$ (i.e. different annihilators of simple
submodules of $M$ are incomparable), we have $Ann(M)\neq \bigcap_{K\in
B\backslash \{N\}}Ann(K)$ for every $N\in B$. The result follows now from
Theorem \ref{Theorem 5.17}.$\blacksquare$
\end{Beweis}

\begin{ex}
\label{Example 5.20}Consider $M=\mathbb{Z}_{30}[x]$ as a $\mathbb{Z}$%
-module. Let $K_{i}=(10x^{i})$, $N_{i}=(15x^{i})$ and $L_{i}=(6x^{i})$. Set $%
K:=\bigoplus\limits_{i=1}^{\infty }K_{i}$, $N:=\bigoplus\limits_{i=1}^{%
\infty }N_{i}$ and $L:=\bigoplus\limits_{i=1}^{\infty }L_{i}$. Notice that
\begin{equation*}
M=K+N+L.
\end{equation*}%
It is clear that $M$ is second representable semisimple with infinite
length, and
\begin{equation*}
att^{s}(M)=Min(att^{s}(M))=\{(2),(3),(5)\}.
\end{equation*}%
Since $K$ is second but not simple, $_{\mathbb{Z}}M$ is not comultiplication
by Theorem \ref{Theorem 5.17} (notice also that $_{\mathbb{Z}}M$ is not
multiplication).
\end{ex}

\begin{ex}
Consider $M=\mathbb{Z}_{30}=(10)+(6)+(15)$. It is clear that $M$ is a second
representable, multiplication, comultiplication and semisimple $\mathbb{Z}$%
-module in which $att^{s}(M)=Min(att^{s}(M))$ and every second submodule of $%
M$ is simple. By Corollary \ref{Corollary 5.19}, every PS-hollow submodule
of $M$ is simple, and so $(10),(6)$ and $(15)$ are the only PS-hollow
submodules of $M$.
\end{ex}

\begin{thm}
\label{Theorem 5.21}

\begin{enumerate}
\item If $M=\sum\limits_{i=1}^{n}K_{i}$ is a minimal second representation
of $M$ with $att^{s}(M)=Min(att^{s}(M))$ and $K_{i}\cap \sum\limits_{j\neq
i}K_{j}$ is PS-hollow in $M$ for all $i\in \{1,\cdots ,n\}$, then $%
M=\bigoplus\limits_{i=1}^{n}K_{i}$ if and only if $K_{i}\cap K_{j}=0$ for
all $i\neq j$.

\item Let $_{R}M$ be distributive and $M=\sum\limits_{i=1}^{n}K_{i}$ be a
minimal PS-hollow representation such that every submodule of $K_{i}$ is
zero or strongly irreducible or $H_{i}$-PS-hollow. Then $M=%
\bigoplus_{i=1}^{n}K_{i}$.
\end{enumerate}
\end{thm}

\begin{Beweis}
\begin{enumerate}
\item Assume that $K_{i}\cap K_{j}=0$ for all $i\neq j$ in $\in \{1,\cdots
,n\}$. Set $I_{i}=\bigcap\limits_{j\neq i}Ann(K_{i})$. Since $%
att^{s}(M)=Min(att^{s}(M))$, we have $I_{i}M=K_{i}$. Also, $K_{i}\cap
\sum\limits_{j\neq i}K_{j}\subseteq K_{i}$. Since $K_{i}\cap
\sum\limits_{j\neq i}K_{j}$ is PS-hollow and each $K_{j}=I_{j}M$ for all $%
j\neq i$, we have $K_{i}\cap \sum\limits_{j\neq i}K_{j}\subseteq
\sum\limits_{j\neq i}K_{j}$ implies that $K_{i}\cap \sum\limits_{j\neq
i}K_{j}\subseteq K_{l}$ for some $l\neq i$, whence $K_{i}\cap
\sum\limits_{j\neq i}K_{j}\subseteq K_{l}\cap K_{i}=0$.

\item Since $_{R}M$ is distributive, it is enough to prove that $K_{i}\cap
K_{j}=0$ for all $i\neq j$ in $\{1,\cdots ,n\}$. Suppose that $K_{i}\cap
K_{j}\neq 0$ for some $i\neq j$. But $0\neq K_{i}\cap K_{j}\subseteq K_{i},$
whence $K_{i}\cap K_{j}$ is strongly irreducible or $H_{i}$-PS-hollow.
Suppose that $K_{i}\cap K_{j}$ is strongly irreducible. Since $K_{i}\cap
K_{j}\subseteq K_{i}\cap K_{j}$, it follows that $K_{i}\subseteq K_{i}\cap
K_{j}$ or $K_{j}\subseteq K_{i}\cap K_{j}$ and so $K_{i}\subseteq K_{j}$ or $%
K_{j}\subseteq K_{i}$ which contradicts the minimality of $%
\sum_{i=1}^{n}K_{i}$. So, $K_{i}\cap K_{j}$ is $H_{i}$-PS-hollow and at the
same time $H_{j}$-PS-hollow, which contradicts the minimality of $%
\sum\limits_{i=1}^{n}K_{i}$. Therefore $K_{i}\cap K_{j}=0$ for all $i\neq j$
in $\{1,\cdots ,n\}.$$\blacksquare$
\end{enumerate}
\end{Beweis}

\begin{exs}
\label{Example 5.22}

\begin{enumerate}
\item Every second representable semisimple module satisfies the assumptions
of Theorem \ref{Theorem 5.21} ( 2).

\item $M=\mathbb{Z}_{n},$ considered as a $\mathbb{Z}$-module, $M$ satisfies
all assumptions of Theorem \ref{Theorem 5.21} ((1) and (2)).
\end{enumerate}
\end{exs}

\begin{thm}
\label{Theorem 5.23}Let $R$ be Artinian and $M=\sum\limits_{i=1}^{n}K_{i}$
be a minimal PS-hollow representation of $_{R}M$. Suppose that the
submodules of $K_{i}$ are PS-hollow $\forall i\in \{1,\cdots ,n\}$. If $%
In(K_{i})\cap In(K_{j})=0$ $\forall i\neq j$ in $\{1,\cdots ,n\},$ then $%
M=\bigoplus\limits_{i=1}^{n}K_{i}$.
\end{thm}

\begin{Beweis}
Assume that $In(K_{i})\cap In(K_{j})=0$ for all $i\neq j$ in $\{1,\cdots
,n\}.$ For each $j\in \{1,\cdots ,n\}$, set $N_{j}:=K_{j}\cap \sum_{i\neq
j}K_{i}$. Then $N_{j}\subseteq In(K_{i})$ for some $i\neq j$. Otherwise, $%
N_{j}\nsubseteq In(K_{i})$ for all $i\neq j$, and so for all $i\neq j$ there
is $I_{i}\in H_{i}$ such that $N_{j}\nsubseteq I_{i}M$. But $N_{j}\subseteq
\sum\limits_{i\neq j}K_{i}\subseteq I_{i}M$ and $N_{j}$ is a PS-hollow
submodule by assumption, whence $N_{j}\subseteq I_{i}M$ for some $i\neq j$
in $\{1,\cdots ,n\}$ (a contradiction).

Observe that $N_{j}\subseteq K_{j}\subseteq In(K_{j})$ and so $%
N_{j}\subseteq In(K_{i})\cap In(K_{j})$ for some $i\neq j$ in $\{1,\cdots
,n\}.$ It follows that $N_{j}=0$ for all $j\in \{1,\cdots ,n\}$ and
therefore $M=\bigoplus\limits_{i=1}^{n}K_{i}$.$\blacksquare$
\end{Beweis}

\begin{cor}
\label{Corollary 5.24}Let $R$ be Artinian and $M=\sum\limits_{i=1}^{n}K_{i}$
a minimal PS-hollow representation of $_{R}M$. Suppose that the nonzero
submodules of $In(K_{i})$ are $H_{i}$-PS-hollow for all $i\in \{1,\cdots
,n\}$, where $K_{i}$ is $H_{i}$-PS-hollow for each $i\in \{1,\cdots ,n\}$.
Then $M=\bigoplus\limits_{i=1}^{n}K_{i}$.
\end{cor}

\begin{Beweis}
Suppose that $In(K_{i})\cap In(K_{j})\neq 0$ for some $i\neq j$ in $%
\{1,\cdots ,n\}.$ Then $In(K_{i})\cap In(K_{j})$ is $H_{i}$-PS-hollow, and
at the same time $In(K_{i})\cap In(K_{j})$ is $H_{j}$-PS-hollow, which is a
contradiction since $H_{i}\neq H_{j}$ as $M=\sum\limits_{i=1}^{n}K_{i}$ is a
minimal PS-hollow representation. Therefore $In(K_{i})\cap In(K_{j})=0$. The
result is obtained by Theorem \ref{Theorem 5.23}.$\blacksquare$
\end{Beweis}

\addcontentsline{toc}{section}{\protect\numberline{}{Index}} \printindex

\end{document}